\documentclass[11pt]{article}
\usepackage{amsmath}
\usepackage{amssymb}
\usepackage{amsthm}
\usepackage{epsf}

\newcounter{AbcT}

\newtheorem {Theorem}    {Theorem}[section]

\newtheorem {Lemma}      [Theorem]    {Lemma}
\newtheorem {Corollary}  [Theorem]    {Corollary}

\newtheorem {Conjecture}[Theorem]     {Conjecture}

\newcommand{\ignore}[1]{}

\def\P{{\bf{P}}}

\def\N{{\bf N}}

\def\C{{\cal{C}}}
\def\K{{\cal{K}}}

\def\RR{\bf R}

\def\|{\, | \, }
\def\v0{{\bf 0}}

\def\0{\hat{0}}
\def\1{\hat{1}}

\def\path{{\tt path}}

\def\phi{\varphi}

\def\be{\begin{equation}}
\def\ee{\end{equation}}

\newcommand{\Z}{\mathbb{Z}}
\newcommand{\halb}{\frac{1}{2}}
\newcommand{\prob}{{\bf P}}
\newcommand{\D}{{\cal D}}
\newcommand{\en}{{\cal E}}

\title{The diameter of long-range percolation clusters on finite cycles
}
\author{Itai Benjamini and Noam Berger\footnote{Research  partially
supported by NSF grant
 \#DMS-9803597 and
by a US-Israel BSF grant. Part of the research was done while
this author was at the Hebrew University of Jerusalem. }}
\begin{document}
\maketitle

\begin{abstract}
Bounds for the diameter and expansion of the graphs created by  long-range
percolation on the cycle $\Z/N\Z$, are given.

\end{abstract}

\section{Introduction}

Model the metric of the world at the year $1900$ by some graph, e.g. $\Z^2$
or $\Z^3$,
equipped with the graph metric. The introduction of fast
communication and transportation starting with phones, cars, airplanes and
finally
(for now) the Internet, decreases distances.
This can be modeled by adding new edges between far away vertices in one
way or another. A natural way to do that, is long-range percolation. In
(Bernoulli)
long-range percolation, a countable set of vertices $V$ is given,
equipped with a distance function $d$, on the set. Now to get a random graph
with $V$ as its vertices set, attach  an edge
$e_{v,u}$ between $v, u \in V$ with probability $p_{d(v,u)}$, determined
only by the
distance between $v$ and $u$, independently of all other pair of vertices.
Long-range percolation on $\Z$ was introduced and studied in \cite{s},
\cite{NS} and
\cite{AN}. These papers mainly studied when an infinite cluster exists,
whether it is unique (yes \cite{uniq})
and the type of phase transitions that occurs, see also \cite{ms}.
Not much attention was given to the geometry and structure of the infinite
cluster, once it exists. In \cite{Be} and \cite{BB} the random walk on and volume
growth of long-range percolation
clusters on
$\Z$ and $\Z^2$ were
studied.
We then observed that even on finite graphs long-range percolation might
be of
interest. In particular in trying to study the world wide web. The
spatial
structure of the world still manifests itself in the web structure. These
ideas
are not new, see \cite{wa} \cite{sw}, and appear occasionally  under the
name "small world".

The attachment probabilities we will consider will have polynomial decay in
the
distance, i.e. $p_{d(v,u)} \sim \beta d^{-s}$.
In the next section the model is defined, in Section
\ref{sec.diam} we consider the diameter of the clusters, Section \ref{moreon} contains the
formulation
of a sharp result (quoted from \cite{BKPS}) regarding the diameter of the
cluster for
$s < 1$.
In Section \ref{chheg} we discuss expansion properties of the cluster.
We end with some less formal concluding remarks.

Similar models have been discussed in the computer science and physics literature.
In \cite{KB} Kleinberg has
studied the properties of a quite similar model. His interest regarded a
two-dimensional fixed degree model for which the
probability of  a long range edge of length $d$ to be
open is proportional to $d^{-s}$. His focus was on constructing
good routing algorithms that rely only on local information.

In \cite{Jesper}, the behavior of random walk on models of this sort is
studied. Jespersen and Blumen (\cite{Jesper}) study the return probabilities of the random
walker in a slightly different model. Their study, as well as \cite{Be},  reveals
a phase transition at $s=2$. The same phase transition (as well as one at $s=1$) shows also
as a result of our  study. A very interesting continuation of this work
was recently written by Coppersmith, Gamarnik and Sviridenko (\cite{Gamar}).

\section{The model}\label{model}

The model we discuss is the finite long-range percolation model  with
polynomial decay.
Let $N$ be a positive integer, let $s,\beta>0$, and consider the following
random graph:
\\The vertices are the elements of the cycle $\Z/N\Z$. Define
$\rho(x,y)=\min(|x-y|,N- |x-y|)$.
\\Determine the edges in the graph as follows:
\\If $\rho(x,y)=1$, then $x$ and $y$ will be attached to each other.
Otherwise, if $x\neq y$, then $x$ and $y$
will be attached with probability $1-\exp(-\beta\rho(x,y)^{-s})$. The
different edges are all independent of
each other. The probability of an edge between to (distant enough) vertices
is very close to
$\beta\rho(x,y)^{-s}$, and this, as well as independence, are the two
important features of the
presented distribution.

We call the graph created this way $G_{s,\beta}(N)$.

Other interesting models could be high dimensional models,
models with different decay rate (exponential or other), or
models with dependencies, like the Random Cluster Model or models
of the type discussed in \cite{KB}.

\section{The diameter}\label{sec.diam}
In this section, we fix $s$ and $\beta$, and estimate the diameter of the
graph.
{}For $N$, $s$ and $\beta$, let $\D(N)=\D_{s,\beta}(N)$ be the diameter of
$G_{s,\beta}(N)$.
The main results are
\begin{Theorem}\label{diam}
(A) If $s>2$ then there is $C=C_{s,\beta}$ s.t.
$$\lim_{N\to\infty}\prob(\D(N)<CN)=0.$$ Moreover, there exists $0<\eta\leq 1$ s.t.
$\D(N)/N\to\eta$ in distribution.
 \\(B) If $s<1$ then there is constant $C=C_{s,\beta}$ s.t.
$$\lim_{N\to\infty}\prob(\D(N)>C)=0$$.
 \\(C) If $1<s<2$ then there is constant $\delta=\delta_{s,\beta}$ s.t.
$$\lim_{N\to\infty}\prob(\D(N)>\log(N)^\delta)=0$$.
\\(D) if $1<s<2$, then there is a constant $C=C_{s,\beta}$ s.t.
$$\lim_{N\to\infty}\prob(\D(N)<C\log(N))=0$$.
\end{Theorem}


At a previous version of this paper, we have conjectured:
\begin{Conjecture}\label{conj_false}
(A) If $s=2$, then the diameter's order of magnitude is $N^\delta$, where
$\delta$ is a function of $\beta$.
 \\(B) If $s=1$, then the diameter's order of magnitude is $\log(N)$.
 \\(C) If $1<s<2$ then the diameter is $\theta(\log(N)^\gamma )$ where
$\gamma>1$ is a function of $s$.
\end{Conjecture}

Recently, Coppersmith, Gamarnik and Sviridenco (\cite{Gamar}) have proved that for $s=1$,
the diameter is of order of magnitude $\log(N)/\log\log(N)$. This contradicts part (A)
of conjecture \ref{conj_false}. At the same paper
they also proved part (B) of Conjecture \ref{conj_false} for $\beta<1$. For all other values
of $\beta$, they have proved an upper bound for the diameter of order of magnitude
$N^\delta$ where $\delta<1$ and its value depends on $\beta$. Further, Biskup
(\cite{mbiskup}) has recently announced that when $1<s<2$ the diameter
is $\log(N)^{\log_2(2/s)+o(1)}$. This proves part (C) of conjecture
\ref{conj_false}.

In view of these results, we now believe:
\begin{Conjecture}\label{conj_true}
(A) If $s=2$, then the diameter's order of magnitude is $N^\delta$, where
$\delta$ is a function of $\beta$.
 \\(B) If $1<s<2$ then the diameter is $\theta(\log(N)^\gamma )$ where
$\gamma>1$ is a function of $s$.
\end{Conjecture}

In \cite{BB} it is shown that for $s=2$ and any $\beta$ the diameter is no less than
$N^{1/\log\log N}$.

\begin{proof}[Proof of Theorem \ref{diam}]
(A) Assume first that the model is the line $[0,N-1]$ with $\rho(x,y)= |x-y|$ and not the circle.
Now, for given $x_0$, the probability that there is no edge between any $x<
x_0$ and any $y>x_0$
is
\begin{eqnarray*}
\prod_{0\leq x< x_0,N>y>x_0}{e^{-\beta(y-x)^{-s}}}
  &=& e^{-\beta\sum_{0\leq x< x_0,N>y>x_0}{(y-x)^{-s}} }          \\
  &\geq& e^{-\beta\sum_{x< x_0,y>x_0}{(y-x)^{-s}} }          \\
  & =  & e^{-\beta\sum_{k=1}^\infty{k^{1-s}} }      \\
  & =  & \psi(s,\beta) > 0.    \\
\end{eqnarray*}
Define a {\em cut} to be a vertex with this property. Take $C$ s.t.
$6C=\halb\psi(s,\beta)$.
So, by the ergodic theorem, if $N$ is large enough then, with
probability as high as we like, there are at least
$6CN = \halb\psi(s,\beta)$ cuts, and therefore the diameter is, with the
same probability, at least $6CN$. Returning  to the cycle, we can divide it
into two lines, of length
$\halb N$ each. Each of these halves is of diameter at least $3CN$. The same
kind of
calculation yields that, with high probability, the
edges between the two halves of the cycle don't reach the middle third of
each of the lines of the vertices, and therefore the diameter stays above
$CN$. In order to prove that the limit of $\D(N)/N$ exists, we do the following:
First, consider long-range percolation on $\Z$. $\D'(N)$ will be the diameter of
the long-range percolation restricted to $[0,N]$. By the sub-additive ergodic theorem,
$D'(N)/N\to\eta$ a.s. for some $\eta>0$.
In order to prove the convergence for the diameter of the cycle, we divide the cycle
$\Z/N\Z$ into two intervals $I_1$ and $I_2$ of length $N/2$. The diameter of each half is
(with very high probability) approximately $\eta N/2$. The longest connection between the
two halves is of length $o(N)$, so there are cut points $x_1,x_2\in I_1$ and $x_3,x_4\in I_2$
s.t. $\rho(x_1,x_2)=N/2-o(N)$ and $\rho(x_3,x_4)=N/2-o(N)$, there are no edges between the
arcs $[x_1,x_2]$ and $[x_3,x_4]$, and, by the strong Markov property, the diameter of each
of these arcs is $\eta N/2-o(N)$. So, we are done. It is of interest to study the fluctuations
of $D(N)/N$ from $\eta$.

(B) The graph dominates the $G(n,p)$ random graph with edge probability
$\beta N^{-s}$. It is known
(see, e.g. \cite{JLR}) that there exists a constant $C$ s.t.
\begin{eqnarray}\label{bd}
\lim_{n\to\infty}\prob(\D(G(n,\beta n^{-s}))<C)=1
\end{eqnarray}
Since the diameter is a decreasing function (w.r.t the standard partial
order), (\ref{bd}) applies
also for our model with $s<1$.

Actually, in this case we can even say more: The infinite graph whose
vertices are the integers,
s.t. every two vertices are attached with probability $1-\exp(-\beta
G|x-y|^s)$ has, a.s.,
a finite diameter, see next section.

(C) Here we use an argument in the spirit of Newman and Schulman's
renormalization
(see \cite{NS}):
Again, assume that the model is a line instead of a circle. This assumption
creates a measure
which is dominated by the original one, and therefore it suffices to prove
the result for the line.
Take
$$C_i=e^{\alpha^i},$$
where $\alpha>1$ is s.t.
$$s<(2-s)\sum_{i=1}^{\infty}\alpha^{-i}$$
Let $k_0$ be a large number, and define $\gamma$ to be
$$\gamma=(2-s)\sum_{i=1}^{k_0}\alpha^{-i}-s.$$

Taking $\alpha$ small enough and $k_0$ large enough,  we can get that
\begin{eqnarray}\label{gamalp}
e^\gamma >\alpha
\end{eqnarray}

Now, take
$$N_k=\prod_{i=1}^{k}C_k.$$

We divide the interval of length $N_k$ into $C_k$ intervals of length
$N_{k-1}$. Each of
these, we divide into $C_{k-1}$ intervals of length $C_{k-2}$, and so on.
This structure has a lot in common
with the one used in \cite{NS} for proving the existence of the
infinite cluster.
We use the following terminology: The $N_l/N_k$ intervals of
length $N_k$ obtained
by this division from $N_l$, are called {\em components of degree $k$}. The
components of
degree $k-1$ inside such a component are {\em sub-components}.
two intervals (or component) $I$ and $J$ are said to be
{\em attached} to each other, if there exists a bond between a point in $I$
and a point in
$J$.

It is enough to show that for some constant $D$,
\begin{equation}\label{nk}
\lim_{k\to\infty}\prob(\D(N_k)\leq D^k)=1,
\end{equation}
because
$$
D^k=\left(\alpha^k\right)^\delta<\log(N_k)^\delta
$$
for $\delta=\frac{\log D}{\log\alpha}$, and
if $N_k<n<N_{k+1}$, then we can bound
$\D(n)$ by $\D(N_{k+1})$. This will
be enough because for some constant $\kappa$, we
have $\log(N_k)\leq\kappa\log(n)$.

We now prove (\ref{nk}): Take some $\epsilon>0$. We will show that for $l$
large enough,
$
\prob(\D(N_l)\leq D^l)>1-\epsilon
$.
Define $\beta_k=\halb\beta N_k^{2-s}$.  Then the probability that two
intervals of length $N_k$ of distance
$lN_k$ from each other have an edge between them is at least
$1-\exp(-\beta_kl^{-s})$.

Take $k_1>k_0$ so large that for every $k>k_1$, we have $\frac{1}{6}\beta
e^{\gamma\alpha^k}>2\alpha^k$, and so that
$e^{-k_1}<\epsilon/2$.
{}For every $k>k_1$, consider the line of length $N_k$. Divide it into $C_k$
components of size
$N_{k-1}$. The probability that not all of the $C_k$ components are attached
to each other is
bounded by
\begin{eqnarray*}
{\binom{C_k}{2}}e^{\left(-\beta_{k-1}C_k^{-s}\right)}
&\leq& C_k^2e^{-\beta_{k-1}C_k^{-s}}\\
&=& C_k^2\exp\left(-\halb\beta
e^{-s\alpha^k}e^{(2-s)\sum_{i=1}^{k-1}\alpha^i}\right)\\
&\leq& \exp\left(2\alpha^k-\halb\beta e^{\gamma\alpha^k}\right) \\
&\leq& \exp\left(-\frac{1}{3}\beta e^{\gamma\alpha^k}\right)
\end{eqnarray*}

Now, take $l$  s.t. $\psi\log l>k_1$ where $\psi$ is s.t.
$\psi\log\alpha=1$. Consider the following event,
denoted by
$\nu$: for every  $\psi\log l\leq k\leq l$, and for every component of degree $k$,
all of its sub-components of
degree $k-1$ are attached to each other.

Given $\nu$, the diameter is no more than $$2^lN_{\psi\log
l}<2^l(e^{(1/(1-\alpha^{-1}))\alpha^{\psi\log l}})<D^l$$
for $D=2\exp(1/(1-\alpha^{-1}))$. Therefore, we want to estimate the
probability of $\nu$:

Take $k$ s.t. $\psi\log l\leq k\leq l$. The probability that there exist a
component of degree $k$, s.t. not all
of its sub-components of degree $k-1$ are attached to each other could be
bounded by the
number of components of degree $k$ times the probability for this event at
each of them, i.e. by
\begin{eqnarray*}
N_k\cdot\exp\left(-\frac{1}{3}\beta e^{\gamma\alpha^k}\right)&\leq&N_l\cdot\exp\left(-\frac{1}{3}\beta e^{\gamma\alpha^k}\right)\\
&\leq&\exp\left(
(1/(1-\alpha^{-1}))\alpha^l-\frac{1}{3}\beta e^{\gamma\alpha^{\psi\log l}}
\right)\\
&=&\exp\left(
(1/(1-\alpha^{-1}))\alpha^l-\frac{1}{3}\beta e^{\gamma l}
\right)\\
&\leq&\exp\left(-\frac{1}{4}\beta e^{\gamma l}
\right)
\end{eqnarray*}
by (\ref{gamalp}) for large enough $l$.

Therefore, the probability of $\nu$ is at least
$$1-l\exp\left(-\frac{1}{4}\beta e^{\gamma l}\right)>1-\epsilon$$
for $l$ large enough.

So,
$$
\D(N)=O\left(\log(N)^
{\frac{\log(2\exp((1-\frac{2s}{2+s})^{-1}))}{\log(2+s)-\log(2s)}+\epsilon}\right)
$$
for every $\epsilon$.


(D) When $s>1$, the expected value of the number of vertices attached to a
certain vertex is
finite. Therefore, the graph's growth rate is bounded by the growth
rate of a Galton-Watson
tree. Thus, its growth rate
is (bounded by) exponential, and so the diameter cannot be smaller than a
logarithm of the number
of vertices ($N$, in this case).
\end{proof}

\section{More on $s < 1$}\label{moreon}

In this section we will report on  a theorem from \cite{BKPS}  for
long-range
percolation, that deals with the case $s <1$ .

Denote by $B^d_N$ the standard $d$-dimensional square lattice restricted
to the $N \times N....\times N \mbox{ }$ box. The diameter of $B^d_N$,
in the graph metric,
denoted $d(\mbox{ } ,\mbox{ } )$, is $dN$. We will add random edges to
$B^d_N$
as follows.
{}Fix some $s \in (0,d)$, and for any two vertices $v, u \in B^d_N$,
add an edge between $v,u$ with probability $d(v,u)^{-s}$, independently
from all other edges. Denote by $D(B^d_N)$ the random diameter
of the box once the new edges were added.

\begin{Theorem}\label{s<1}(\cite{BKPS})

$$
\lim_{N \to \infty}D(B^d_N) =
\lceil \frac{d}{ d-s}\rceil \mbox{ }  a.s.
$$
\end{Theorem}

\noindent{\bf Remarks}.
\begin{itemize}
\item
This theorem was proved in \cite{BKPS} as a corollary of the general theory of
stochastic dimension that was developed for the study of uniform spanning forest.
A random relation  $\RR \subset \Z^d \times \Z^d$ is said to have
{\em stochastic dimension} $s$, if there is some constant $c>0$ such
that
for all $x\neq y$ in $\Z^d$,
$$
c^{-1}|x-y|^{s-d}<\,\P[x \RR y]<c |x-y|^{s-d},
$$
and certain correlation inequalities
hold.
The results regarding stochastic dimension are formulated and proven
in this generality, to allow application in several contexts. One application
is the above. Another application is the following:

{}For every $v\in B^d_N$, let $S^v$ be a simple random walk
starting from $v$, with $\{S^v\}_{v\in B^d_N}$ independent.
Let $v \mbox{ knows } u$ be the relation $\exists n\in\N\ \ S^v(n)=S^u(n)$.
Then $\mbox{``knows''}$ has stochastic dimension $2$. The "know" diameter (i.e. the diameter
of the graph in which there is an edge between $v$ and $u$ whenever $v$ knows $u$) of
$B^d_N$ will be $\lceil \frac{d}{d-2}\rceil \mbox{ }  a.s.$

\item
If indeed we are all six handshakes
away from any other person on the planet, then assuming $d=2$, the
"real world handshakes exponent" $s$, might be around $1- 2/7$.

\end{itemize}

\section{Cheeger's Constant}\label{chheg}
We would like to explore the Cheeger constant of these graphs. First,
recall its
definition. For a set of vertices $A$ in a finite graph, let
$\partial A$ be the set of edges
$\{e = (v,u), v \in A, u \in A^c\}$.
A geometric tool which is used in order to bound relaxation times  of
random walk on graphs is the Cheeger constant $\C(G)$, see \cite{AF}.
\begin{equation} \label{eq:def_phi}
\C(G) = \inf_A \frac{\partial A}{|A|},
\end{equation}
where the infimum is taken over all non empty set of vertices $A$, with
$|A| \leq \frac{|G|}{2}$.

Again, fix $s$ and $\beta$, and denote
by $\C(N)=\C_{s,\beta}(N)$ the Cheeger constant of the graph
$G_{s,\beta}(N)$
when $s>2$, the Cheeger constant is $\theta(\frac{1}{N})$. That is proven the same
way as the fact that the diameter is
linear in $N$. However, for $1<s\leq2$, the Cheeger constant exhibits an
interesting behavior:
\begin{Theorem}\label{cheeger}
(A) If $1<s<2$, then there is a constant $\alpha=\alpha(s)$ s.t.
$$
\lim_{N\to\infty}\prob\left(\C(N)>N^{-\alpha}\right)=0.
$$
(B) If $s=2$, then for every $\epsilon$, there $C$ s.t.
$$
\limsup_{N\to\infty}\prob\left(\C(N)>\frac{C\log N}{N}\right)<\epsilon.
$$
\end{Theorem}
Part (B) is proved as lemma 3.4 of \cite{Be}.
\begin{proof}[Proof of (A)]
Divide the circle into two arcs, $A$ and $B$,  of length $\halb N$ each. The
expected value of the
number of edges connecting the two halfs is:
\begin{eqnarray*}
2&+&\sum_{i\in A,j\in B}(1-\exp(-\beta\rho(i,j)^{-s}))\\
&\leq&2+2\sum_{k=2}^{\halb N}k(1-\exp(-\beta k^{-s}))\\
&\leq&2+C\sum_{k=2}^{\halb
N}k^{1-s}\leq2+2C\beta\sum_{k=2}^{\infty}k^{1-s}\\
&\leq&2+\frac{C}{2-s}\left(\frac{N}{2}+1\right)^{2-s}
\end{eqnarray*}
{}For some constant $C$.
So, the expected value of the size of the boundary of $A$ divided by the
size of $A$ is bounded by
$
DN^{1-s}
$
for some constant $D$, and this, using the Markov inequality,  gives the
desired result for
every $\alpha<s-1$.
\end{proof}

As a simple corollary of Theorem \ref{cheeger}, we get the following lower bound
for the mixing time of a random walk on $G_{s,\beta}(N)$. (See \cite{AF} for
background on mixing)

\begin{Corollary}\label{mixing}
Denote by $\tau(G)$ the mixing time of the simple random walk on a graph
$G$. If $1<s<2$ then for every $\delta<s$, we have
$$
\lim_{n\to\infty}\prob(\tau(G_{s,\beta}(N))<N^\delta)=0
$$
\end{Corollary}

By theorem  \ref{cheeger} we see that the size of the smallest cut in the
case
$1<s<2$ reminds that of a cube in $Z^d$.  However, the cut structure is
different.
While in a cube of length $N^{1/d}$ (and,
therefore, $N$ vertices) in $Z^d$ one can find
a sequence of length $\theta(N^{1/d})$ of nested cuts of size
$O(N^\frac{d-1}{d})$ each, in $G_{s,\beta}(N)$ the length of a sequence of nested cuts is bounded by the
diameter which is no more than poly-logarithmic in $N$.

\section{Concluding Remarks}\label{conc}

\subsection{Discussion}
The geometry of the $1\leq s<2$ clusters described here is different from the
geometry of other natural graphs - In the $s > 1$ case
its diameter is rather short
(poly-logarithmic in
the volume),
while its smallest cut is also small - as small as that of a box in
the
$n$-dimensional lattice (Cut sets  polynomial in the volume).

For example, when $s=1$, the average degree is $\theta(\log(n))$, as was shown in \cite{Gamar},
the diameter is $\theta\left(\frac{\log(n)}{\log\log(n)}\right)$. This might lead to thinking
that this graph is similar to the random graph $G(n,\log(n)/n)$. However, there exist
large sets (such as the vertices $[1,...,\halb n]$) that have small boundaries
of order $n/logn$. This can be seen by following the argument in the proof of Theorem
\label{cheeger} part (A).

\subsection{electrical resistance}
By the proof of Lemma 2.4 of  \cite{Be}, if $1<s<2$, there is a constant $C$ which does not depend
on $N$, s.t. if we pick at random two vertices of $G_{s,\beta}(N)$, then with a very high probability the
effective electrical resistance between them is bounded by $C$. However, the maximal electrical
resistance is unbounded - in fact it is easy to see that it is at least logarithmic with $N$.

\subsection{Inverse problems}

As we saw above, once the tail of the connecting probabilities
is fat (but not too fat) we get a graph of poly-logarithmic  diameter and
super-polynomial volume
growth. This brings the question whether the geometry of the underlying
graph, say $\Z^d$ for different values of $d$, is disappearing  and we get
some
universal generic geometries?
To be more precise, assume you are given a sample from a super critical
long-range percolation taking place
on one of the two graphs $\Z$ or $\Z^2$, with some connecting probabilities,
which are not given
to you and without the labeling of the edges by $\Z$ or $\Z^2$.
Can you a.s. tell if the sample came from long-range percolation on $\Z$?
I.e. is the set of measures on graphs coming from considering all super
critical
long-range percolations
on $\Z$ is singular with respect to that coming from $\Z^2$?
The question could be asked for any pair of graphs, in particular can one
distinguish between $\Z$ and the $3$-regular tree $T_3$, or between
or between $\Z^{d_1}$ and $\Z^{d_2}$?  Another variant: Replace a graph by a
graph roughly isometric -- can one distinguish $\Z$ from $\Z^2$ but not from
the
triangular lattice on $\Z^2$?

More generally, let $G$ be a graph, and assume you are given a sample from a
super critical long-range percolation
on $G$, with some unknown connecting probabilities and without the labeling
of the vertices of the sample by their names in $G$. What information can
be recovered on $G$? For instance can properties such as  amenability of
$G$,
transitivity of $G$, number of ends in $G$, volume growth of $G$,
planarity of $G$, and isoperimetric dimension of $G$,
be recovered?

\noindent
Itai Benjamini,\\
The Weizmann Institute and Microsoft research\\
e-mail:itai@wisdom.weizmann.ac.il\\
\\
Noam Berger, \\
The University of California\\
e-mail:noam@stat.berkeley.edu\\
\end{document}